\documentclass[11pt,a4paper]{article}
\usepackage[dvips]{graphicx}
\usepackage{color}
\usepackage[dvips]{epsfig}
\usepackage{latexsym}
\usepackage{makeidx}
\newcommand{\qed}{\hfill \mbox{$\Box$} \vspace{.5\baselineskip}}

\newcommand{\RR}{{\sf I\hspace{-0,08em} R}}                     

\newcommand{\parx}{\frac{\partial}{\partial x}}

\newcommand{\G}{\mathcal{G}}      
\newcommand{\Gcon}{\mathcal{C}}      

\newcommand{\B}{\mathcal{B}}    
\newcommand{\Bprime}{\mathcal{B}^{\prime}}    

\newcommand{\CB}{C_{\mathcal{B}}}       
\newcommand{\CBprime}{C_{\mathcal{B}}^{\prime}}       
\newcommand{\CBpoint}{C_{\mathcal{B}}^{\bullet}}       
\newcommand{\CBB}{C_{\mathcal{B}}^{\diamondsuit}}  
\newcommand{\CBSB}{C_{\mathcal{B}}^{\bullet\negthinspace\diamondsuit}}

\newcommand{\ar}{\mbox{\Large{$a$}}}  

\newcommand{\Hn}{\mathop{\rm H}\nolimits}

\newcommand{\Hu}{\mathop{\rm Hu}\nolimits}      
\newcommand{\Hup}{\Hu^{\bullet}}      

\newcommand{\Ca}{\mathop{\rm Ca}\nolimits}  

\newcommand{\Oc}{\mathop{\rm Oc}\nolimits}
\newcommand{\Ocp}{\mathop{\rm Oc}^{\bullet}\nolimits}

\newcommand{\yg}{\textbf{y}}

\newcommand{\KK}{{\sf I\hspace{-0,08em} K}}                     

\usepackage{color}

\def\yes{yes}

\def\incol{\color{red}}
\def\outcol{\color{blue}}

\def\UseMapleColor#1{%
\ifx#1\yes\gdef\incol{\color{red}}\gdef\outcol{\color{blue}}%
\else\edef\incol{}\edef\outcol{}%
\fi}

\newlength{\abovemapleskip}
\newlength{\belowmapleskip}
\newenvironment{maplelike}{%
\par%
\vskip\abovemapleskip%
\bgroup%
\parindent=0mm%
\parskip=0mm%
\abovedisplayskip=0mm%
\abovedisplayshortskip=0mm%
\belowdisplayskip=0mm%
\belowdisplayshortskip=0mm%
}{%
\vskip\belowmapleskip%
\egroup%
}

\newlength{\tmplgtha}
\newcommand{\mwsin}[1]{%
\par
\bgroup%
\incol
\hskip2mm{\footnotesize$>$}\hskip3mm%
\settowidth{\tmplgtha}{\hskip2mm{\footnotesize$>$}\hskip3mm}%
\tmplgtha=-\tmplgtha%
\advance\tmplgtha by\textwidth%
\begin{minipage}[t]{\tmplgtha}%
{\tt#1}%
\end{minipage}%
\egroup%
\par%
}

\newcommand{\mwsout}[1]{%
\outcol#1%
}%

\newlength{\mapleind}
\newtheorem{theorem} {Theorem} [section]
\newtheorem{proposition} [theorem] {Proposition}

\title{Enumerative problems inspired by \\Mayer's theory of cluster integrals}
\author{Pierre Leroux, LaCIM, UQAM}  
\date{December 31, 2003} 
\begin{document}

\maketitle

\begin{abstract} 
The basic functional equations for connected and 2-connnected graphs can be traced
back to the statistical physicists Mayer and Husimi.
They play an essential role in establishing rigorously the virial expansion
for imperfect gases. We survey this approach and inspired by these equations, we investigate 
the problem of enumerating some classes of connected graphs all of whose blocks are contained 
in a given class $B$. 
Included are the species of Husimi graphs ($B =$ "complete graphs"), 
cacti ($B =$ "unoriented cycles"), and oriented cacti ($B =$ "oriented cycles"). 
For each of these, we consider the question of their labelled 
or unlabelled enumeration and of their molecular expansion, 
according (or not) to their block-size distributions.
\end{abstract}
\begin{resume}
Les \'equations fonctionnelles de base pour les graphes connexes et 2-connexes
remontent aux physiciens Mayer et Husimi. Ces relations sont importantes
pour \'etablir rigoureusement le d\'eveloppement du viriel pour les gaz imparfaits.
Nous pr\'esentons cette d\'emarche et inspir\'es par ces relations, 
nous examinons le probl\`eme du d\'enombrement de quelques classes de graphes
connexes dont les blocs sont pris dans une famille $B$ donn\'ee. 
Cela inclut les graphes de Husimi ($B =$ "graphes complets"), 
les cactus ($B =$ "polygones"), et les cactus orient\'es ($B =$ "cycles orient\'es"). 
Pour chacune de ces esp\`eces, on s'int\'eresse au d\'enombrement \'etiquet\'e, 
ordinaire ou selon la distribution des tailles des blocs,
au d\'enombrement non-\'etiquet\'e, et au d\'eveloppement mol\'eculaire.
\end{resume}%
%
%
\section{Introduction}

\subsection{Combinatorial species and functional equations for connected graphs and blocks}

Informally, a \emph{combinatorial species} is a class of labelled discrete
structures which is closed under isomorphisms induced by relabelling along
bijections. See Joyal \cite{Jo81} and Bergeron, Labelle, Leroux \cite{BLL98} for an
introduction to the theory of species. 
To each species $F$ are associated a number of series expansions among which are 
the (exponential) generating function, $F(x)$, for labelled enumeration,
defined by 
\begin{equation} \label{eq:expgf}
F(x) = \sum_{n\geq 0} |F[n]| \frac {x^n}{n!},
\end{equation}
where $|F[n]|$ denotes the number of $F$-structures on the set $[n]=\{1,2,\ldots,n\}$,
the (ordinary) \emph{tilde} generating function $\widetilde{F}(x)$, for unlabelled enumeration, 
defined by
\begin{equation} \label{eq:ordgf}
\widetilde{F}(x) = \sum_{n\geq 0} \widetilde{F}_n {x^n},
\end{equation}
where $\widetilde{F}_n$ denotes the number of isomorphism classes of F-structures
of size $n$, %
the \emph{cycle index series} $Z_F(x_1,x_2,x_3,\cdots)$, defined by 
\begin{equation} \label{eq:zf}
Z_F(x_1,x_2,x_3,\cdots) = \sum_{n\geq 0} \frac{1}{n!} \left(\sum_{\sigma\in S_n}
\mathop{\rm fix} F[\sigma]\ x_1^{\sigma_1}x_2^{\sigma_2}x_3^{\sigma_3}\cdots \right),
\end{equation}
where $S_n$ denotes the group of permutations of $[n]$, $\mathop{\rm fix} F[\sigma]$
is the number of $F$-structures on $[n]$ left fixed by $\sigma$,
and $\sigma_j$ is the number of cycles of length $j$ in $\sigma$,
and the molecular expansion of $F$, which is a description of the $F$-structures 
and a classification according to their stabilizers and will be discussed later.

Combinatorial operations are defined on species: sum, product, (partitional) composition, 
derivation, 
which correspond to the usual operations on the exponential generating functions. 
And there are rules for computing the other associated series, involving plethysm.
See \cite{BLL98} for more details.
An equality $F=G$ between species is a family of bijections between structures,
$$\alpha_U : F[U] \rightarrow G[U],$$
where $U$ ranges over all underlying (labelling) sets,
which commutes under relabellings. 
Such an identity gives rise to an equality between all the series expansion
associated to species.
\begin{figure}[ht] 
\begin{center}
\includegraphics[scale=.9]{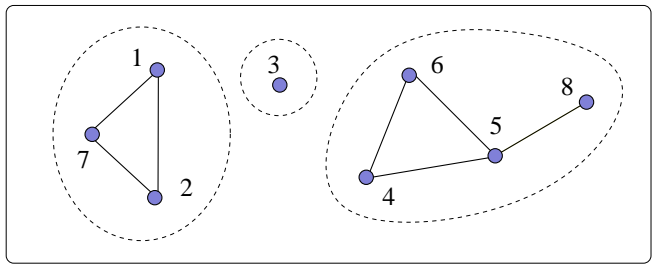}
\end{center}
\caption{A simple graph $g$ and its connected components \label{fig:multcompconnexe}}
\end{figure}

For example, the fact that any simple graph on a set (of vertices) $U$ is the
disjoint union of connected simple graphs (see Figure \ref{fig:multcompconnexe}) 
is expressed by the equation
\begin{equation} \label{eq:graph}
\G=E(\Gcon),
\end{equation}
where $\G$ denotes the species of (simple) graphs, $\Gcon$, that of connected graphs,
and $E$, the species of \emph{Sets} (in French: \emph{Ensembles}).
There correspond the well-known relations for their exponential generating functions,
\begin{equation} \label{eq:egfgraphe}
\G(x)= \exp(\Gcon(x))
\end{equation}
and for their tilde generating functions,
\begin{eqnarray} 
\widetilde{\G}(x) & = & Z_E(\widetilde{\Gcon}(x),\widetilde{\Gcon}(x^2),\ldots) \nonumber\\
    & = & \exp{\left(\sum_{k\geq1} \frac{1}{k} \widetilde{\Gcon}(x^k)\right)}.
               \label{eq:tildegraphe}
\end{eqnarray}

\begin{figure}[ht] 
\begin{center}
\includegraphics[scale=.8]{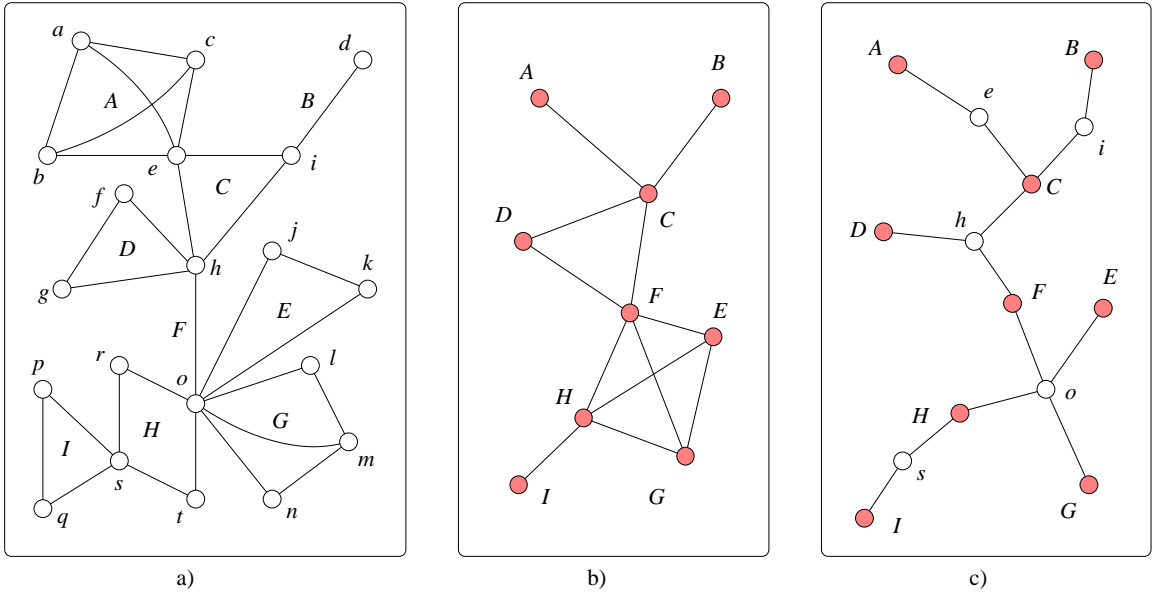}
\end{center}
\caption{a) A connected graph $g$, b) the block-graph of $g$, 
c) the block-cutpoint tree of $g$  \label{fig:graphblock}}
\end{figure}

\noindent
{\bf Definitions.} A \emph{cutpoint} (or \emph{articulation point}) of a
connected graph $g$  is  a vertex of $g$ whose removal yields a disconnected 
graph. A connected graph is called \emph{2-connected} if it has no cutpoint.
A \emph{block} in a simple graph is a maximal 2-connected subgraph.
The \emph{block-graph} of a graph $g$ is a new graph whose vertices are
the blocks of $g$ and whose edges correspond to blocks having a common
cutpoint. The \emph{block-cutpoint tree} of a connected graph $g$ is a graph 
whose vertices
are the blocks and the cutpoints of $g$ and whose edges correspond to incidence
relations in $g$. See Figure \ref{fig:graphblock}.
%

Now let $\B$
be a given species of 2-connected graphs. We denote by  $\CB$ the species of
connected graphs all of whose blocks are in $\B$, called $\CB$-graphs. 

\vspace{1mm}
\noindent {\bf Examples 1.1.}
Here are some examples for various choices 
of $\B$:
\begin{enumerate}

\item 
If $\B = \B_a$, the class of \emph{all} 2-connected graphs, then $\CB=\Gcon$, 
the species of (all) connected graphs.
\item 
If  $\B=K_2$, the class of "edges", then $\CB=\ar$, the species of (unrooted, free) trees
($\ar$ for French \emph{arbres}).
\item 
If $\B=\{P_m, m \geq 2\}$, where $P_m$ denotes the class of size-$m$ polygons (by convention,
$P_2=K_2$), then $\CB=\Ca$, the species of cacti. A \emph{cactus} can also be defined as a  
connected graph in which no edge lies in more than one cycle. 
Figure \ref{fig:cactus}, a), represents a typical cactus.
\begin{figure}[bht] 
\begin{center}
\includegraphics[scale=.8]{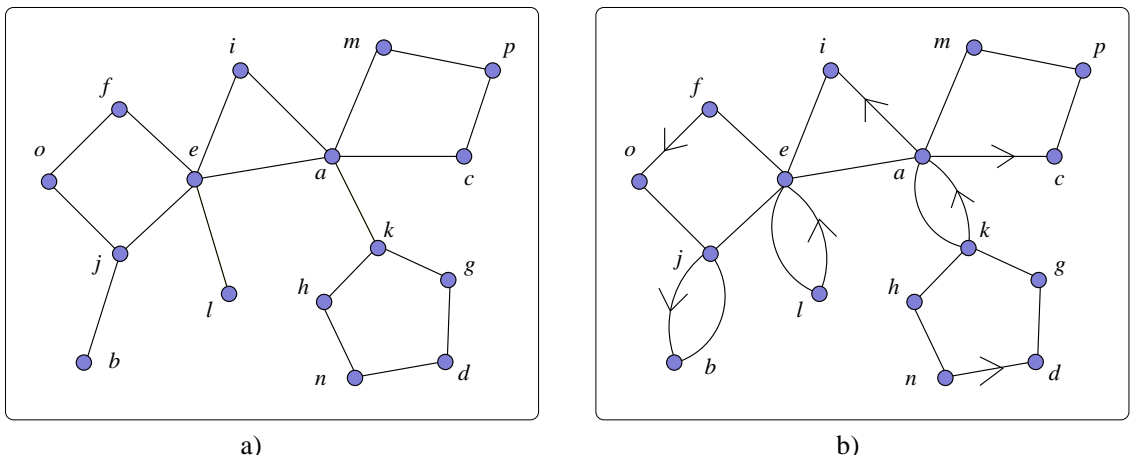}
\end{center}
\caption{a) a
typical cactus, b) a typical oriented cactus \label{fig:cactus}}
\end{figure}
\item 
If $\B=K_3=P_3$, the class of "triangles", then $\CB=\delta$, the class of triangular cacti.
\item  
If $\B=\{ K_n, n \geq 2 \}$, the family of complete graphs, then $\CB= \Hu$, the species
of \emph{Husimi graphs}, that is of connected graphs whose blocks are complete graphs. They were
first (informally) introduced by Husimi in \cite{Hu50}.
A Husimi graph is shown in Figure \ref{fig:graphblock}, b).
See also Figure \ref{fig:husimigraph}. It can be easily shown that any Husimi graph is 
the block-graph of some connected graph.
\item 
If $\B = \{ C_n, n \geq 2 \}$, the family of oriented cycles, then  $\CB = \Oc$, the species of 
\emph{oriented cacti}. Figure \ref{fig:cactus}, b) shows a typical oriented cactus. 
These structures were introduced by C. Springer \cite{Sp96} in 1996.
Although directed graphs are involved here, the functional equations (\ref{eq:cbprime}) 
and (\ref{eq:cbpointw}) given below are still valid. 

\end{enumerate}

\noindent
{\bf Remark.} Cacti where first called Husimi trees. See for example \cite{HaNo53}, \cite{HaUl53},
\cite{Ri51} and \cite{UhFo63}. However this term received much criticism since they are not 
necessarily trees. Also, a careful reading of Husimi's article \cite{Hu50} shows that the graphs
he has in mind and that he enumerates (see formula (\ref{eq:hun2n3bis}) below) are the Husimi graphs defined
in  item 5 above. The term cactus is now widely used. See Harary and Palmer \cite{HaPa73}.
Cacti appear regularly in the mathematical litterature, for example in
the classification of base matroids \cite{MaZa03}, and in combinatorial optimization 
\cite{Fl99}. 

The following functional equation (see (\ref{eq:cbprime})) is fairly well known. It can be found in various forms
and with varying degrees of generality in \cite{BLL98}, \cite{HaPa73}, \cite{Jo81},
\cite{La83}, \cite{Le88}, \cite{LeMi92}, \cite{No54}, \cite{Ri51}, \cite{Ro70}.
In fact, it was anticipated by the physicists (see \cite{Hu50} and \cite{UhFo63})
in the context of Mayers' theory of cluster integrals as we will see below.
The form given here, in the structural language of species, is the most general one since all the
series expansions follow. It is also the easiest form to prove.

Recall that for any species $F=F(X)$, the derivative $F^{\prime}$ of $F$ is the species defined
as follows: an $F^{\prime}$-structure on a set $U$ is an $F$-structure on the set $U\cup\{\ast\}$,
where $\ast$ is an external (unlabelled) element. In other words, one sets 
$$F^{\prime}[U] = F[U + \{\ast\}].$$
Moreover, the operation $F\mapsto F^{\bullet}$, of pointing (or rooting) $F$-structures 
at an element of the underlying set, can be defined by 
$$ F^{\bullet} = X\cdot F^{\prime}.$$

\begin{theorem}
Let $\B$ be a class of 2-connected graphs and $\CB$, the species
of connected graphs all of whose blocks are in $\B$. We then have the functional equation
\begin{equation} \label{eq:cbprime}
\CBprime=E(\Bprime(\CBpoint)).
\end{equation}
\end{theorem}
\begin{figure}[ht] 
\begin{center}
\includegraphics[scale=.6]{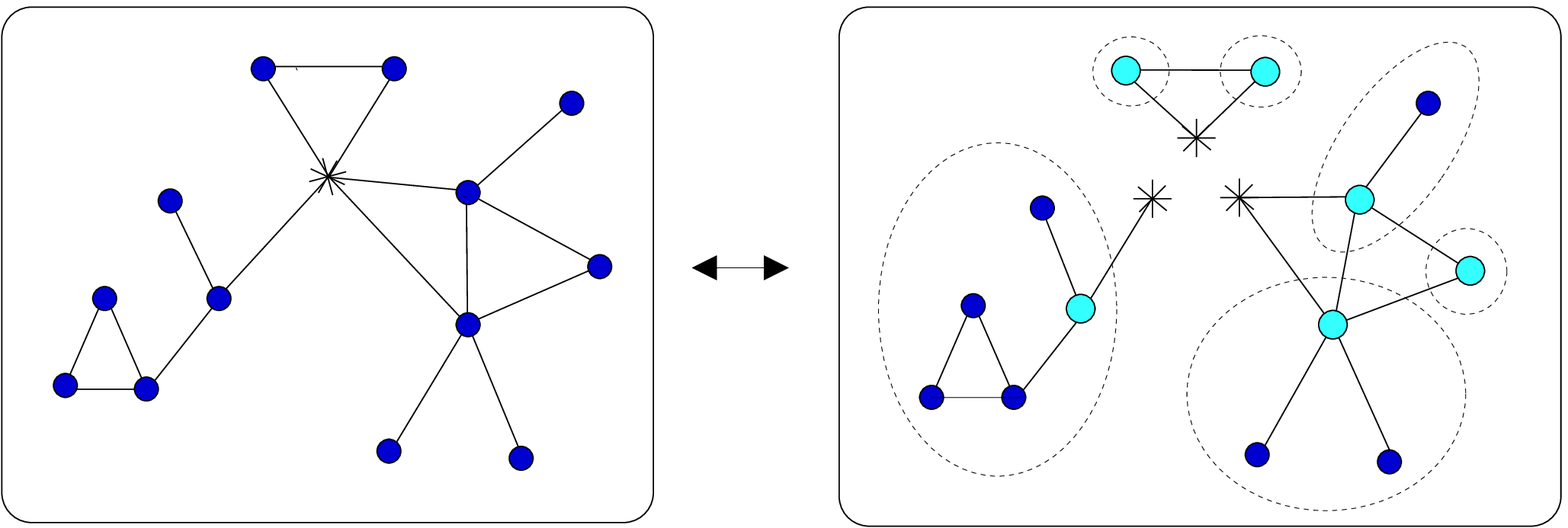}  
\end{center}
\caption{$\CBprime=E(\Bprime(\CBpoint))$ \label{fig:fig424} }
\end{figure}

\noindent {\bf Proof.} See Figure \ref{fig:fig424}. \qed

Multiplying (\ref{eq:cbprime}) by $X$, one finds 
\begin{equation} \label{eq:cbpoint}
\CBpoint = X \cdot E(\Bprime(\CBpoint)),
\end{equation}
and, for the exponential generating function,
\begin{equation} \label{eq:cbpointexp}
\CBpoint(x) = x \cdot \exp(\Bprime(\CBpoint(x))).
\end{equation}
%
\subsection{Weighted versions}

Weighted versions of these equations are needed in the applications. See for example
Uhlenbeck and Ford \cite{UhFo63}. A \emph{weighted species} is a species $F$ 
together with weight functions $w=w_U:F[U]\rightarrow \KK$ defined on $F$-structures,
which commute with the relabellings.  
Here $\KK$ is a commutative
ring in which the weights are taken, usually a ring of polynomials or formal power series
over a field of characteristic zero. 
We write $F=F_w$ to emphasize the fact that $F$ is a weighted species 
with weight function $w$.
The associated generating functions are
then adapted by replacing set cardinalities $|A|$ by total weights 
$$|A|_w = \sum_{a\in A}w(a).$$
The basic operations on species are also adapted to the weighted context,
using the concept of Cartesian product of weighted sets: 
Let $(A,u)$ and $(B,v)$ be weighted sets. A weight function $w$ is defined
on the Cartesian product $A\times B$ by 
$$w(a,b) = u(a)\cdot v(b).$$
We then have \ \  $|A\times B|_w = |A|_u \cdot |B|_v$.

\vspace{.5\baselineskip}
\noindent {\bf Definition.}
A weight function $w$ on the species $\G$ of graphs is said to be \emph{multiplicative
on the connected components} if for any graph $g\in\G[U]$, whose connected components are
 $c_1,c_2,\ldots,c_k$, we have
$$w(g)=w(c_1)w(c_2)\cdots w(c_k).$$

\vspace{.5\baselineskip}
\noindent{\bf Examples 1.2.} 
The following weight functions $w$ on the species of graphs are multiplicative
on the connected components.
\begin{enumerate}
  \item $w_1(g) := y^{e(g)}$, where 
   $e(g)$  is the number of edges of $g$. 
  \item $w_2(g)$  =  graph complexity of $g$ :=  number of maximal spanning forests of $g$. 
  \item $w_3(g) := x_0^{n_0}x_1^{n_1}x_2^{n_2}\cdots$, where $n_i$ is the number of 
vertices of degree $i$. 
\end{enumerate}
%

\begin{theorem}
Let $w$ be a weight function on graphs which is multiplicative on the connected components.
Then we have 
\begin{equation}
\G_w=E\left(\Gcon_w\right).  \label{eq:graphw}
\end{equation}
\end{theorem}
For the exponential generating functions, we have
$$G_w(x)=\exp(\Gcon_w(x)),$$
where
$G_w(x) \ = \ \sum_{n\geq 0}|\G[n]|_w \frac{x^n}{n!}
  \ = \ \sum_{n\geq 0} (\sum_{g\in G[n]}w(g)) \frac{x^n}{n!},$
and similarly for $\Gcon_w(x).$

\vspace{.5\baselineskip}
\noindent {\bf Definition.}
A weight function on connected graphs is said to be \emph{block-multiplicative} 
if for any connected graph $c$, whose blocks are $b_1,b_2,\ldots,b_k$, we have
$$w(c)=w(b_1)w(b_2)\cdots w(b_k).$$ 
\noindent{\bf Examples 1.3.} 
The weight functions
$w_1(g)=y^{e(g)}$\ and \  $w_2(g) =$ graph complexity of $g$
of Examples 1.2 are block-multiplicative,
but not the function
$w_3(g)=x_0^{n_0}x_1^{n_1}x_2^{n_2}\cdots$.
Another example of a block-multiplicative weight function is obtained by 
introducing formal variables $y_i$ ($i\geq2$) marking the block sizes.
In other terms, if the connected graph $c$ has $n_i$ blocks of size $i$, 
for $i=2,3,\ldots$, one sets 
\begin{equation}
w(c)=y_2^{n_2}y_3^{n_3}\cdots. \label{eq:yiweight}
\end{equation}

The following result is then simply the weighted version of Theorem 1.1.

\vspace{.5\baselineskip}
\noindent {\bf Theorem 1.3.} 
Let $w$ be a block-multiplicative weight function on
connected graphs whose blocks are in a given species $B$.  Then we have
\begin{equation} \label{eq:cbpointw}
\left({\CBpoint}\right)_w = X\cdot E(\Bprime_{w}((\CBpoint)_w)).
\end{equation}

\subsection{Outline}

In the next section, we see how equations (\ref{eq:graphw}) and (\ref{eq:cbpointw})
are involved in the thermodynamical study of imperfect (or non ideal) gases, following Mayers'
theory of cluster integrals \cite{Mama40}, as presented in Uhlenbeck and Ford 
\cite{UhFo63}.  In particular, the virial expansion, which is a kind of asymptotic
refinement of the perfect gases law, is established rigourously, at least in its
formal power series form. See equation (\ref{eq:virialexp}) below. 
It is amazing to realize that the coefficients of the virial expansion involve 
directly the total valuation $|\B_a[n]|_w$, for $n\geq 2$, of 2-connected graphs. 
An important role in this theory is also played by the
enumerative formula (\ref{eq:hun2n3bis}) for labelled Husimi trees according to their block-size
distribution, which extends Cayley's formula $n^{n-2}$ for the number of labelled trees of size $n$.

Motivated by this, we first consider, in Section 3, the labelled 
enumeration of some classes of connected graphs of the form $\CB$, according
or not to their block-size distribution.  Included are the species of Husimi graphs, cacti,
and oriented cacti. The methods involve the Lagrange inversion formula and
Pr\"ufer-type bijections.
It is also natural to examine the question of unlabelled enumeration of these structures. 
This is a more difficult problem, for two reasons.  First, equation (\ref{eq:cbpointw})
deals with rooted structures and it is necessary to introduce a tool for counting
the unrooted ones. Traditionally, this is done by extending the Dissimilarity charactistic
formula for trees of Otter \cite{Ot48}. See for example \cite{HaNo53}.
Inspired by formulas of Norman (\cite{FoNoUh56}, (18)) and Robinson (\cite{Ro70}, Theorem 7),
we have given over the years a more 
structural formula which we call the Dissymmetry theorem for graphs, whose proof is
remarkably simple and which can esily be adapted to various classes of tree-like structures; 
see \cite{BLL98}, \cite{BBLL00}, \cite{FGLL02}, \cite{La92}, \cite{LaLaLe02}, 
\cite{LaLaLe03}, \cite{Le88}, \cite{LeMi92}. Second, as for trees, it should not be expected
to obtain simple closed expressions but rather recurrence formulas for the number
of unlabelled $\CB$-structures. Three examples are given in that section.

Finally, in Section 4, we present the molecular expansion of some
of these species. This expansion can be computed first for the rooted species, 
by iteration, and the dissymmetry theorem is then invoked for the unrooted ones.
The computations can be carried out easily using the Maple package "Devmol"
developped at LaCIM and available at the URL www.lacim.uqam.ca; see \cite{AuLe03}.

\vspace{.5\baselineskip}
\noindent {\bf Acknowledgements.} 
This talk is partly taken from my student M\'elanie Nadeau's "M\'emoire de ma\^{\i}trise" \cite{Na02}. 
I would like to thank her and Pierre Auger for their considerable help,
and also Abdelmalek Abdesselam, Andr\'e Joyal, Gilbert Labelle, Bob Robinson, and Alan Sokal, 
for useful discussions.

\section{
Some statistical mechanics}

\subsection{Partition functions for the non-ideal gas}

\begin{figure}[bht] 
\begin{center}
\includegraphics[scale=.7]{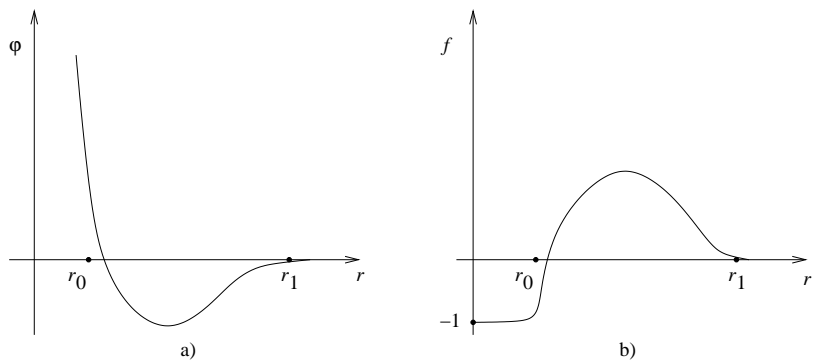}  
\end{center}
\caption{a) the function $\varphi(r)$, b) the function $f(r)$ \label{fig:phietf} }
\end{figure}

Consider a non-ideal gas, formed of $N$ particles interacting in a vessel 
$V\subseteq \RR^3$ (whose volume is also denoted by $V$) and whose positions are 
$\overrightarrow{x_1}, \overrightarrow{x_2}, \ldots,\overrightarrow{x_N}$.
The Hamiltonian of the system is of the form
\begin{equation} 
H=\sum_{i=1}^{N} \left(\frac{\overrightarrow{p_i}^2}{2m}+ U(\overrightarrow{x_i})\right) + 
\sum_{1\leq i<j\leq N} \varphi(|\overrightarrow{x_i}-\overrightarrow{x_j}|), 
\label{eq:hamiltonien}
\end{equation}
where 
$\overrightarrow{p_i}$
is the linear momentum vector and
$\frac{\overrightarrow{p_i}^2}{2m}$ is the kinetic energy of the  $i^{th}$ particle,
$U(\overrightarrow{x_i})$ is the potential at position $\overrightarrow{x_i}$
due to outside forces (e.g., walls),
$|\overrightarrow{x_i}-\overrightarrow{x_j}|=r_{ij}$ 
is the distance between the particles $\overrightarrow{x_i}$ and $\overrightarrow{x_j}$,
and it is assumed that the particles interact only pairwise through the central potential
$\varphi(r)$. 
This potential function $\varphi$ has a typical form shown in Figure 
\ref{fig:phietf} a). 

The \emph{canonical partition function} $Z(V,N,T)$ is defined by
\begin{equation}
Z(V,N,T)=\frac{1}{N!h^{3N}}\int \exp{(-\beta H)} d\Gamma,  \label{eq:canonpartition1}
\end{equation} 
where $h$ is Planck's constant,  $\beta={1\over {kT}}$, $T$ is the absolute temperature
and $k$ is Boltzmann's constant, and $\Gamma$ represents the state space
 $\overrightarrow{x_1},\ldots,\overrightarrow{x_N},\overrightarrow{p_1},
\ldots,\overrightarrow{p_N}$ of dimension $6N$. 
A first simplification comes from the assumption that the potential energy 
$U(\overrightarrow{x_i})$ is negligible or null. Secondly, the integral
over the momenta $\overrightarrow{p_i}$ in (\ref{eq:canonpartition1}) is a product of 
Gaussian integrals which are easily evaluated so that the canonical partition function 
can now be written as
\begin{equation}
Z(V,N,T)=\frac{1}{N!\lambda^{3N}}\int_{V}\cdots\int_{V}
\exp{\left(-\beta \sum_{i<j} \varphi(|\overrightarrow{x_i}-\overrightarrow{x_j}|) \right)} 
d\overrightarrow{x_1}\cdots d\overrightarrow{x_N},
\label{eq:partcanonique2}
\end{equation}
where $\lambda = h(2\pi mkT)^{-\frac{1}{2}}$.

Mathematically, the \emph{grand-canonical distribution} is simply the generating function
for the canonical partition functions, defined by
\begin{equation}
Z_{\mathop{\rm gr}}(V,T,z)=\sum_{N=0}^{\infty} Z(V,N,T)(\lambda^3z)^N,\label{eq:zgrand1}
\end{equation}
where the variable $z$ is called the \emph{fugacity} or the \emph{activity}.
All the macroscopic parameters
of the system are then defined in terms of this grand canonical ensemble. 
For example, the \emph{pressure}, $P$, the average number of particles, $\overline{N}$, 
and the \emph{density}, $\rho$, are defined by
\begin{equation}
\frac{P}{kT}=\frac{1}{V} \log{Z_{\mathop{\rm gr}}(V,T,z)},\ \ \ \label{eq:PsurkT}
\overline{N} = z\frac{\partial}{\partial z} \log{Z_{\mathop{\rm gr}}(V,T,z)},\ \ \ 
\mathrm{and}\ \ \ \rho:=\frac{\overline{N}}{V}.\label{eq:rhoz}
\end{equation}
%
\subsection{The virial expansion}

In order to better explain the thermodynamic behaviour of non ideal gases,  Kamerlingh Onnes
proposed, in 1901, a series expansion of the form
\begin{equation}
\frac{P}{kT} = \frac{\overline{N}}{V} + \gamma_2(T) \left(\frac{\overline{N}}{V}\right)^2
 +\gamma_3(T)\left(\frac{\overline{N}}{V}\right)^3+\cdots, \label{eq:virial}
\end{equation}
called the \emph{virial expansion}. Here $\gamma_2(T)$ is the second virial coefficient,
$\gamma_3(T)$ the third, etc. This expansion was first derived theoretically from the
partition function $Z_{\mathop{\rm gr}}$ by 
Mayer \cite{Mama40} around 1930.
It is the starting point of Mayer's theory of "cluster integrals". 
Mayer's idea consists in setting
\begin{equation}
1 + f_{ij}
=\exp{\left(-\beta \varphi(|\overrightarrow{x_i}-\overrightarrow{x_j}|)\right)},
\label{eq:fij}
\end{equation}
where $f_{ij} = f(r_{ij})$. 
The general form of the function $f(r) = \exp(-\beta\varphi(r))-1$ is shown in Figure
\ref{fig:phietf}, b). 
In particular, $f(r)$ vanishes when $r$ is greater than the range $r_1$ of the 
interaction potential.
By substituting in the canonical partition function (\ref{eq:partcanonique2}), 
one obtains
\begin{equation}
Z(V,N,T)=\frac{1}{N!\lambda^{3N}}\int_{V}\cdots\int_{V} \prod_{1\leq i<j\leq N}
\left(1+f_{ij}\right)d\overrightarrow{x_1}\cdots d\overrightarrow{x_N}.  
\label{eq:partcanonique3}
\end{equation}

The terms obtained by expanding the product 
$\prod_{1\leq i<j\leq N} (1+f_{ij})$ can be represented by simple graphs
where the vertices are the particles and the edges are the chosen factors
$f_{ij}$. The partition function (\ref{eq:partcanonique3}) can then 
be rewritten in the form
\begin{eqnarray}
Z(V,N,T)&=&\frac{1}{N!\lambda^{3N}} \sum_{g\in\G[N]}
   \int_{V}\cdots\int_{V} \prod_{\{i,j\}\in g} f_{ij} \ d\overrightarrow{x_1}
\cdots d\overrightarrow{x_N}  \nonumber \\ 
 & = & \frac{1}{N!\lambda^{3N}} \sum_{g\in\G[N]} W(g), \label{eq:partcanonique4}
\end{eqnarray}
where the weight $W(g)$ of a graph $g$ is given by the integral
\begin{equation} 
W(g)=\int_{V}\cdots\int_{V} \prod_{\{i,j\}\in g} f_{ij} d\overrightarrow{x_1}
\cdots d\overrightarrow{x_N}.  \label{eq:weightofg}  
\end{equation}

For the grand canonical function, we then have
\begin{eqnarray}
 Z_{\mathop{\rm gr}}(V,T,z) & = &\sum_{N=0}^{\infty} Z(V,N,T)(\lambda^3z)^N \nonumber\\
 & = & \sum_{N=0}^{\infty} \frac{1}{N!\lambda^{3N}} 
                    \sum_{g\in\G[N]} W(g)(\lambda^3z)^N \nonumber\\
 & = & \sum_{N=0}^{\infty} \frac{1}{N!} \sum_{g\in\G[N]} W(g)z^N \nonumber\\
 & = & \G_W(z). \label{eq:zgrand2} 
\end{eqnarray}
%
\begin{proposition}
The weight function $W$ is multiplicative on the
connected components.
\end{proposition}

For example, for the graph $g$ of Figure \ref{fig:multcompconnexe},
we have 
\begin{eqnarray*}
\lefteqn{W(g)   =  \int_{V^{8}} f_{12}f_{17}f_{27}f_{45}f_{46}f_{56}f_{58}
             d\overrightarrow{x_1}\cdots d\overrightarrow{x_8} } \\
 & = &\int f_{12}f_{17}f_{27}d\overrightarrow{x_1} 
       d\overrightarrow{x_2} d\overrightarrow{x_7} \int d\overrightarrow{x_3}
       \int f_{45}f_{46}f_{56}f_{58}d\overrightarrow{x_4} d\overrightarrow{x_5} 
        d\overrightarrow{x_6} d\overrightarrow{x_8} \\
 & = & W(c_1)W(c_2)W(c_3), 
\end{eqnarray*}
where $c_1$, $c_2$ and $c_3$ represent the three connected components of $g$.
Following Theorem 1.2, we deduce that
\begin{equation}
\G_W(z)=\exp{\left(\Gcon_W(z)\right)}, \label{eq:gfgraphW}
\end{equation}
where $\Gcon_W$ denotes the weighted species of connected graphs, with
$$
\Gcon_W(z)=\sum_{n\geq1} |\Gcon[n]|_W \frac{z^n}{n!}, 
$$ 
and
\begin{equation}
|\Gcon[n]|_W = \sum_{c\in\Gcon[n]} \int_V\cdots\int_V \prod_{\{i,j\}\in c} f_{ij}
                    \  d\overrightarrow{x_1}\cdots d\overrightarrow{x_n}.  \label{eq:totalweightc} 
\end{equation}
Historically, the quantities $b_n(V) = \frac{1}{Vn!}|\Gcon[n]|_W$
%
are precisely the \emph{cluster integrals} of Mayer. Equation (\ref{eq:gfgraphW}) then 
provides a combinatorial interpretation for the quantity 
$\frac{P}{kT}$. Indeed, one has, by (\ref{eq:PsurkT}),
\begin{eqnarray}
\frac{P}{kT} & = & \frac{1}{V}\log{Z_{\mathop{\rm gr}}(V,T,z)} \nonumber  \\
 & = & \frac{1}{V} \log{\G_W(z)}  \nonumber \\
 & = & \frac{1}{V} \Gcon_W(z). \label{eq:PsurkTbis}
\end{eqnarray}

\begin{proposition}
For $V$ large, the weight function $w(c)=\frac{1}{V}W(c)$,
defined on the species of connected graphs, is block-multiplicative.
\end{proposition}

\noindent {\bf Proof.} First observe that for any connected graph $c$
on the set of vertices $[k]=\{1,2,\ldots\,k\}$, the value of the
partial integral
\begin{equation}
I = I(\overrightarrow{x_k}) = \lim_{V\rightarrow\infty}\,
\int_{V}\cdots\int_{V}\prod_{\{i,j\}\in c}f_{ij} d\overrightarrow{x_1} 
            \cdots d\overrightarrow{x_{k-1}} \label{eq:partialint}
\end{equation}
is in fact independent of $\overrightarrow{x_k}$. Indeed, since the 
$f_{ij}$'s only depend on the relative positions 
$r_{ij} = |\overrightarrow{x_i}-\overrightarrow{x_j}|,$
and considering the short range $r_1$ of the interaction potential
and the connectednes of $c$,
we see that the support of the integrand in (\ref{eq:partialint})
lies in a ball of radius at most $(k-1)r_1$ centered at $\overrightarrow{x_k}$
and that a simple translation
$\overrightarrow{x_i}\mapsto \overrightarrow{x_i}+\overrightarrow{u}$
will give the same value of the integral. It follows that
\begin{eqnarray*}
W(c) & = & \int_{V}\cdots\int_{V}\prod_{\{i,j\}\in c} f_{ij}  
       d\overrightarrow{x_1}\cdots d\overrightarrow{x_{k-1}} \int_{V} d\overrightarrow{x_k} \\
 & \approx & \int_{V}\cdots\int_{V} \prod_{\{i,j\}\in c} f_{ij} 
         d\overrightarrow{x_1}\cdots d\overrightarrow{x_{k-1}}\cdot V, 
\end{eqnarray*}
for $V$ large, and the value of the partial integral (\ref{eq:partialint}) is in fact $w(c)$:
\begin{equation}
\lim_{V\rightarrow\infty}\,\int_V\cdots\int_V  \prod_{\{i,j\}\in c} f_{ij} d\overrightarrow{x_1}\cdots 
d\overrightarrow{x_{k-1}} = \lim_{V\rightarrow\infty}\,\frac{1}{V} W(c) = w(c). \label{eq:weightwc}
\end{equation}

\begin{figure}[bht] 
\begin{center}
\includegraphics[scale=.5]{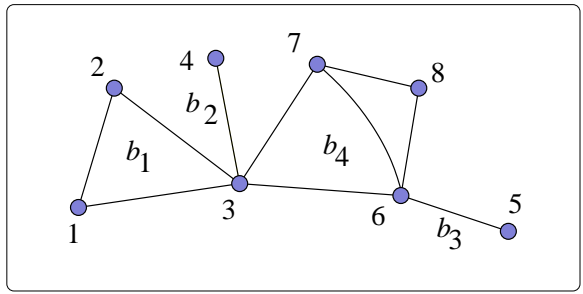}
\end{center}
\caption{A connected graph with blocks $b_1,b_2,b_3,b_4$\label{fig:multcomp2connexe}}
\end{figure}

\begin{figure}[bht] 
\begin{center}
\includegraphics[scale=.5]{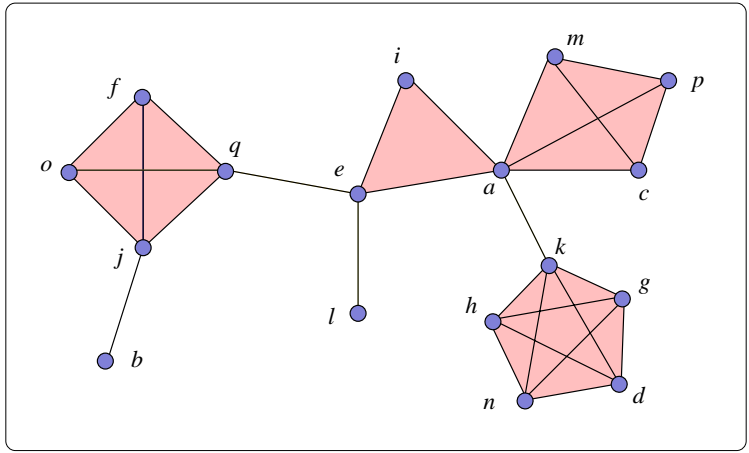} 
\end{center}
\caption{A Husimi graph with block-size distribution $(4,1,2,1,0,\ldots)$
\label{fig:husimigraph}}
\end{figure}

Now if a connected graph is decomposed into blocks $b_1,b_2,\ldots,b_{k}$, we have
$$ 
w(c)=w(b_1)w(b_2)\cdots w(b_k).
$$
For example, for the graph $c$ shown in Figure \ref{fig:multcomp2connexe},
we have
\begin{eqnarray*} 
\lefteqn{w(c)  =  \int_{V^7} f_{12}f_{13}f_{23}f_{34}f_{56}f_{37}f_{36}f_{67}f_{68}f_{78}
   d\overrightarrow{x_1} d\overrightarrow{x_2}\cdots d\overrightarrow{x_7} } \\ 
 & = & \int f_{12}f_{13}f_{23}d\overrightarrow{x_1} d\overrightarrow{x_2} 
       \,f_{34}d\overrightarrow{x_4}\,f_{56}d\overrightarrow{x_5}\,f_{37}f_{36}f_{67}f_{68}f_{78}
       d\overrightarrow{x_3}d\overrightarrow{x_6}d\overrightarrow{x_7} \\
 & = & w(b_1)w(b_2)w(b_3)w(b_4).
\end{eqnarray*} 
\qed

By Theorem 1.3, it follows that
\begin{equation}
\Gcon_w^{\bullet} = X\cdot E(\Bprime_w(\Gcon_w^{\bullet})), \label{eq:cwpoint}
\end{equation}
where  $\B=\B_a$ is the species of all 2-connected graphs, and for the exponential generating
functions,
\begin{equation}
\Gcon_w^{\bullet}(z) = z \exp(\Bprime_w(\Gcon_w^{\bullet}(z))). \label{eq:fgcwpoint}
\end{equation}

\subsection{Computation of the virial expansion}

The virial expansion (\ref{eq:virial}) can now be established, following Uhlenbeck and Ford
\cite{UhFo63}. From (\ref{eq:rhoz}), we have, for the density $\rho(z)=\frac{\overline{N}}{V}$,
\begin{eqnarray}
\rho(z) & = & z\frac{\partial}{\partial z} \frac{1}{V}\log{Z_{\mathop{\rm gr}}(V,T,z)} \nonumber\\
 & = & z\frac{\partial}{\partial z} \Gcon_w(z)  \nonumber\\
 & = & \Gcon_w^{\bullet}(z). \label{eq:rhoeqcpoint}
\end{eqnarray}
Hence $\rho(z)$ satisfies the functional equation (\ref{eq:fgcwpoint}), that is
\begin{equation}
\rho(z)=z\exp{\Bprime_w(\rho(z))}. \label{eq:eqfunrhoz}
\end{equation}
The idea is then to use this relation in order to express $z$ in terms of $\rho$
in $\frac{P}{kT}(z)$, as follows. We have, by (\ref{eq:PsurkTbis}) and (\ref{eq:rhoeqcpoint}),
\begin{eqnarray}
\frac{P}{kT} & = & \frac{1}{V} \log{Z_{\mathop{\rm gr}}(V,T,z)} \nonumber \\
   & = & \Gcon_w(z)  \nonumber \\
   & = & \int_{0}^{z} \Gcon_w^{\prime}(t)\ dt \nonumber \\
   & = & \int_{0}^{z} \frac{\rho(t)}{t}\ dt. \label{eq:intrhoovert}
\end{eqnarray}
Let us make the change of variable 
$$
t=t(r)=r\exp{(-\B_w^{\prime}(r))},
$$
which is the inverse function of $r=\rho(t)$, by (\ref{eq:eqfunrhoz}).
Note that $\rho(0)=0$ and $\rho(z)=\rho$, and also that
$$
dt = [ \exp{(-\B_w^{\prime}(r))}-r\exp{(-\B_w^{\prime}(r))}\cdot\B_w^{\prime\prime}(r)] dr. 
$$
Pursuing the computation of the integral (\ref{eq:intrhoovert}), we have
\begin{eqnarray}
\frac{P}{kT} & = & \int_{0}^{z} \frac{\rho(t)}{t}\ dt \nonumber \\
 & = & \int_{0}^{\rho} \left(1-r\B_w^{\prime\prime}(r)\right)\ dr \nonumber \\
 & = & \rho- \int_{0}^{\rho} r\B_w^{\prime\prime}(r)\ dr \nonumber \\
 & = & \rho- \int_{0}^{\rho} \sum_{n\geq1} n\beta_{n+1} \frac{r^n}{n!}\ dr \nonumber \\
 & = & \rho- \sum_{n\geq2} (n-1) \beta_{n} \frac{{\rho}^n}{n!}, \label{eq:virialexp}
\end{eqnarray}
where we have set $\B_w(r)=\sum_{n\geq2} \beta_n\frac{r^n}{n!}$.
This is precisely the virial expansion, with $\rho=\frac{\overline{N}}{V}$.
Hence the $n^{\mathop{\rm th}}$ virial coefficient, for $n\geq2$,
is given by
\begin{eqnarray}
\gamma_n(T) & = & -\frac{(n-1)}{n!} \beta_n  \nonumber  \\
 & = & -\frac{(n-1)}{n!} |\B[n]|_w. 
\end{eqnarray} 

Mayer's original proof of the virial expansion is more technical, since he is not aware 
of a direct combinatorial proof of equation (\ref{eq:fgcwpoint}). 
The following observation is used: By grouping the connected graphs on $[n]$ whose block
decomposition determines the same Husimi graph on $[n]$, and then collecting all Husimi graphs
having the same block-size distribution, one obtains, using the 2-multiplicativity of $w$,
\begin{eqnarray*}
|C[n]|_w & = & \sum_{c\in\Gcon[n]} w(c)  \\
 & = & \sum_{\{B_1,B_2,\ldots\}\in \Hu[n]} \sum_{\{b_i\in \B[B_i]\}} \prod_{i} w(b_i)  \\ 
 & = & \sum_{\{B_1,B_2,\ldots\}\in \Hu[n]} \prod_{i}  \sum_{b\in \B[B_i]} w(b)  \\
 & = &  \sum_{\{B_1,B_2,\ldots\}\in \Hu[n]}  \prod_{i} \beta_{|B_i|}\\
 & = &  \sum_{{n_2,n_3,\ldots \atop \sum n_i(i-1)=n-1}} 
             \Hu(n_2,n_3,\ldots)\beta_2^{n_2}\beta_3^{n_3}\cdots, 
\end{eqnarray*}
where $\Hu[n]$ denotes the set of Husimi graphs on $[n]$ and
$\Hu(n_2,n_3,\ldots)$ is the number of Husimi graphs on $[n]$ having
$n_i$ blocks of size $i$, for $i\geq 2$.
Mayer then proves the following enumerative formula 
%
\begin{equation}
\Hu(n_2,n_3,\ldots)=\frac{(n-1)!\, n^{\sum n_j-1}}{\prod_{i\geq2}\, (i-1)!^{n_i}\, n_i!},
     \label{eq:hun2n3}
\end{equation}
and goes on proving (\ref{eq:fgcwpoint}) and (\ref{eq:virialexp}) analytically.

Formula (\ref{eq:hun2n3}) is quite remarkable. It is an extension of Cayley's formula
$n^{n-2}$ for the number of trees on $[n]$ (take $n_2=n-1, n_3=0, \ldots$).
It has many different proofs, using, for example, Lagrange inversion  or a Pr\"ufer
correspondence, and gives the motivation for the enumerative problems related to 
Husimi graphs, cacti, and oriented cacti, studied in the next section.

\subsection{Gaussian model}

It is interesting, mathematically, to consider a Gaussian model, where
\begin{equation}
f_{ij}=-\exp{(-\alpha||\overrightarrow{x_i}-\overrightarrow{x_j}||^2)},
\label{eq:gaussian} 
\end{equation}
which corresponds to a soft repulsive potential, at constant temperature. In this case,
all cluster integrals can be explicitly computed (see \cite{UhFo63}): The weight
$w(c)$ of a connected graph $c$, defined by (\ref{eq:weightwc}), has value

\begin{equation}
w(c)=(-1)^{e(c)}\left(\frac{\pi}{\alpha}\right)^{\frac{3}{2}(n-1)}\gamma(c)^{-\frac{3}{2}},
\label{eq:wgaussian} 
\end{equation}
where $e(c)$ is the number of edges of $c$ and $\gamma(c)$ is the graph complexity
of $c$, that is, the number of spanning subtrees of $c$. This formula incorporates
three very descriptive weightings on connected graphs, which are multiplicative
on 2-connected components, namely,
$$w_1(c)=y^{e(c)},\ \  w_2(c)=\gamma(c),$$
already seen in Examples 1.2, and 
$w_3(c)=u^{n-1},$ where $n$ is the number of vertices.

\section{Enumerative results}
In this section, we investigate the enumeration of Husimi graphs, cacti and oriented cacti,
according, or not, to their block size distribution.
Recall that these classes can be viewed as species of connected graphs of the form $\CB$ 
and that the functional equation (\ref{eq:cbpoint}) 
for rooted
$\CB$-structures can be invoked, as well as its weighted version (\ref{eq:cbpointw}).


\subsection{Labelled enumeration}

\begin{proposition} 
The number $\Hn_n=|\Hu[n]|$ of (labelled) Husimi graphs on 
$[n]$, for $n \geq 1$, is given by 
\begin{equation}
\Hn_n = \sum_{k \geq 0} S(n-1,k)\, n^{k-1}, \label{eq:hun}
\end{equation}
where $k$ represents the number of blocks and $S(m,k)$ denotes the Stirling number of the
second kind.
\end{proposition}

\noindent {\bf Proof.} The species $\Hu$ of Husimi graphs is of the form $\CB$ with   
$\B = K_{\geq 2}$, the class of complete graphs of size $\geq 2$. 
From the species point of view, it is equivalent to take $\B = E_{\geq 2}$, 
the species of sets of size $\geq 2$.  We then have
$\B^{\prime} = E_{\geq 1}$, with exponential generating function 
$E_{\geq 1}(x) = {e^x-1}$. Hence the species $\Hup$ of rooted Husimi graphs
satisfies the functional equation

\begin{equation}
\Hup=X \, E(E_{\geq 1}(\Hup)). \label{eq:hupoint}
\end{equation}
which, for the generating function $\Hup(x)$, translates into
\begin{equation}
\Hup(x) = x R(\Hup(x)), \label{eq:hupointx}
\end{equation}
with $R(x) = \exp(e^x-1)$. The Lagrange inversion formula then gives
\begin{eqnarray*}
[x^n]\Hup(x) & = &  \frac{1}{n}[t^{n-1}]\left(e^{e^t-1}\right)^n \\
    & = & \frac{1}{n}[t^{n-1}]e^{n(e^t-1)} \\
    & = & \frac{1}{n}[t^{n-1}]\sum_{k \geq 0} n^k \frac{(e^t-1)^k}{k!} \\
    & = & \frac{1}{n}[t^{n-1}]\sum_{k \geq 0} n^k \sum_{m} S(m,k)\frac{t^m}{m!} \\
    & = & \frac{1}{n}\sum_{k \geq 0} n^k \frac{S(n-1,k)}{(n-1)!} \\
    & = & \sum_{k \geq 0} n^{k} \frac{S(n-1,k)}{n!}.
\end{eqnarray*}

Since we are dealing with exponential generating functions, we should
multiply by $n!$ to get the coefficient of $x^n\over n!$. Moreover
we should divide by $n$ to obtain the number of unrooted Husimi graphs.
In conclusion, we have
$
\Hn_n = \frac{1}{n}n![x^n] \Hup(x) = \sum_{k \geq 0} S(n-1,k) \, n^{k-1}.
$
The fact that $k$ represents the number of blocks will appear more clearly 
in the bijective proof given below. \qed

We now come to Mayer's enumerative formula (\ref{eq:hun2n3}).

\begin{proposition}
\emph{(Mayer \cite{Ma68}, Husimi \cite{Hu50})}
 Let $(n_2,n_3,\ldots)$
be a sequence of non-negative integers and $n=\sum_{i\geq2}n_i(i-1)+1$.
Then the number $\Hu(n_2,n_3,\ldots)$ of Husimi graphs on $[n]$ having
$n_i$ blocks of size $i$ is given by
\begin{equation}
\Hu(n_2,n_3,\ldots)=\frac{(n-1)!}{(1!)^{n_2}n_2!(2!)^{n_3}n_3!\cdots} \, n^{k-1},
\label{eq:hun2n3bis}
\end{equation}
where $k = \sum_{i\geq2} n_i$ is the total number of blocks.
\end{proposition}

\noindent {\bf Proof.} Here we are dealing with the weighted species $\Hu_w$
of Husimi graphs weighted by the function 
$w(h)=y_2^{n_2}y_3^{n_3}\cdots$
describing the block-size distribution of the Husimi graph $h$.
Technically, we should take
$B_w=\sum_{m\geq2}(E_m)_{y_m}$,
where the index $y_m$ indicates that the sets of size $m$ have weight $y_m$,
for which 
$$
B_v(x)=\sum_{m\geq2}y_m{x^m\over{m!}} \ \ \mathop{\rm and}
  \ \ \Bprime_v(x)=\sum_{m\geq1}y_{m+1}{x^m\over{m!}}.
$$
The functional equation (\ref{eq:cbpointw}) then gives 
\begin{equation}
\Hup_w(x) = x\,(\exp\,(\sum_{m\geq1} y_{m+1} {x^m\over{m!}})\circ\Hup_w(x))
\label{eq:hupwx}
\end{equation}
and Lagrange inversion formula can be used. We find 
%
\begin{eqnarray*}
[x^n]\Hup_w(x) & = & \frac{1}{n}[t^{n-1}]\left( \exp(\sum_{m\geq1} y_{m+1} 
              {t^m\over{m!}})\right)^n 
 \\
 & = & \frac{1}{n}[t^{n-1}]\exp(n\sum_{m\geq1} y_{m+1} {t^m\over{m!}})  \\
 & = & \frac{1}{n}[t^{n-1}]\sum_{k \geq 0} \frac{n^k}{k!}\left(\sum_{m \geq 1} 
        y_{m+1}\frac{t^m}{m!}\right)^k \\
 & = & \frac{1}{n}[t^{n-1}]\sum_{k \geq 0} \frac{n^k}{k!}
          \left(y_2\frac{t}{1!}+y_3\frac{t^2}{2!} +y_4\frac{t^3}{3!}\cdots\right)^k  \\
 & = & \frac{1}{n}[t^{n-1}]\sum_{k \geq 0} \frac{n^k}{k!} \sum_{ k_1+k_2+\cdots=k} 
      \frac{k!}{k_1!k_2!\cdots} \left(y_2\frac{t}{1!}\right)^{k_1}
                         \left(y_3\frac{t^2}{2!}\right)^{k_2}\cdots \\
 & = & \frac{1}{n}[t^{n-1}]\sum_{k \geq 0} \, \sum_{k_1+k_2+\cdots=k}
          \frac{n^k}{k_1!k_2!\cdots} \frac{y_2^{k_1}y_3^{k_2}\cdots}{(1!)^{k_1}(2!)^{k_2}\cdots}
                      t^{k_1+2k_2+\cdots}  \\
 & = & \sum_{k \geq 0} \sum_{{k_1+k_2+\cdots =k \atop k_1+2k_2+\cdots =n-1}} 
       \frac{n^{k-1}}{k_1!(1!)^{k_1}k_2!(2!)^{k_2}\cdots} y_2^{k_1}y_3^{k_2}\cdots .
\end{eqnarray*}
Extracting the coefficient of the monomial
$y_2^{n_2}y_3^{n_3}\cdots$ 
and multiplying by  $(n-1)! = \frac{1}{n} n!$ then gives the result. \qed  

There are many other proofs of (\ref{eq:hun2n3bis}). Husimi \cite{Hu50} establishes a
recurrence formula which, in fact, characterizes the numbers $\Hu(n_2,n_3,\ldots)$.
He then goes on to prove the functional equation (\ref{eq:hupwx}).
Mayer gives a direct proof
which becomes more convincing when coupled with a Pr\"ufer-type correspondence.
Such a correspondence is given by Springer in \cite{Sp96} for the number $\Oc(n_2,n_3,\ldots)$
of labelled oriented cacti having block size distribution $2^{n_2}3^{n_3}\cdots$.

It is easy to adapt Springer's bijection to Husimi graphs. To each Husimi graph $h$
on $[n]$ having $k$ blocks is assigned a pair $(\lambda,\pi)$, where $\lambda$ is
a sequence $(j_1,\ldots,j_{k-1})$ of elements of $[n]$ of length $k-1$, and $\pi$
is a partition of the set $[n]\backslash \{j_{k-1}\}$ into $k$ parts. Moreover,
if $h$ has block size distribution $2^{n_2}3^{n_3}\cdots$, then $\pi$ has part-size
distribution $1^{n_2}2^{n_3}\cdots$.  This is done as follows. A \emph{leaf-block} $b$
of a Husimi graph $h$ is a block of $h$ containing exactly one articulation point, 
denoted by $j_b$. Let $b=b(h)$ be the leaf-block of $h$ for which the set 
$b(h)\backslash \{j_{b(h)}\}$ contains the smallest element among all sets of the form 
$b\backslash \{j_b\}$.

\begin{figure}[bht] 
\begin{center}
\includegraphics[width=9cm,height=10cm]{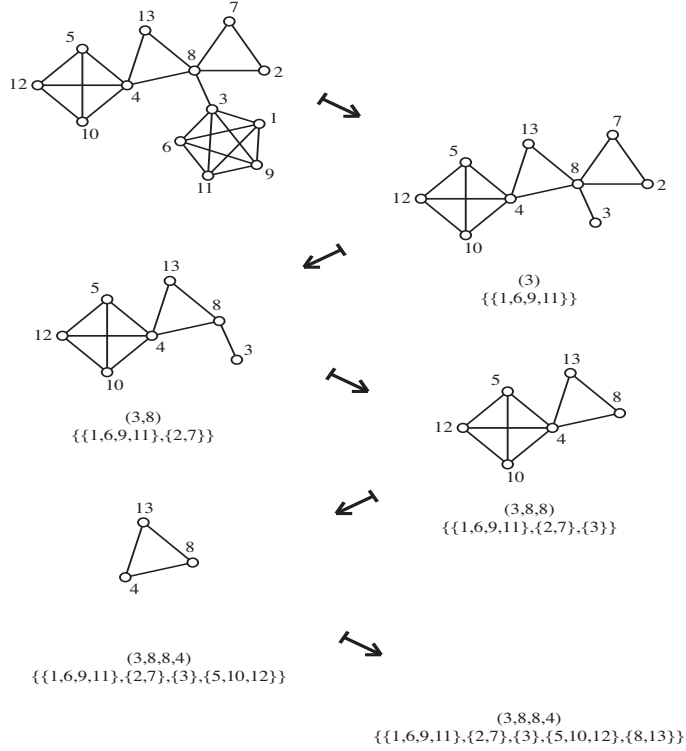}
\end{center}
\caption{Pr\"ufer correspondence for a Husimi graph} \label{fig:bijhypHus}
\end{figure}

 The correspondence proceeds recursively by 
\begin{enumerate}
\item
adding $j_{b(h)}$ to the sequence $\lambda$,
\item
adding the set $b(h)\backslash \{j_{b(h)}\}$ to the partition $\pi$, 
\item
removing the block $b(h)$ (but not the articulation point $j_{b(h)}$) from $h$, and
\item
pursuing with the remaining Husimi graph.
\end{enumerate}
The procedure stops after the $(k-1)^{\mathop{\rm th}}$ iteration when, in supplement, the last
remaining block minus the $(k-1)^{\mathop{\rm th}}$ articulation point $j_{k-1}$ is added to the
partition $\pi$.  An example is shown in Figure \ref{fig:bijhypHus}.
This procedure can easily be reversed and the resulting bijection 
proves both (\ref{eq:hun}) and (\ref{eq:hun2n3bis}).

We now turn to the species of oriented cacti, defined in Example 1.1.6.  These structures
were introduced by Springer in \cite{Sp96} for the purpose of enumerating "ordered short 
factorizations" of a circular permutation $\rho$ of length $n$
into circular permutations, $a_i$ of length $i$ for each $i$. Such a factorization 
is called \emph{short} if $\sum_{i\geq2}(i-1)a_i = n-1$.

\begin{proposition}
The number $\Oc_n=|\Oc[n]|$ of oriented cacti on $[n]$, for $n \geq 2$, is given by

\begin{equation}
\Oc_n= \sum_{k \geq 1} \frac{(n-1)!}{k!} {n-2\choose k-1} \, n^{k-1}, \label{eq:ocn}
\end{equation}
where  $k=\sum_{i\geq 2}n_i$ is the number of cycles.
\end{proposition}

\noindent {\bf Proof.} The proof is similar to that of (\ref{eq:hun}). Here 
$\B=\sum_{k\geq2}C_k$, the species of oriented cycles of size $k\geq2$ 
and $\Bprime = \sum_{k\geq1} X^k$, the species of "lists" (totally ordered sets)
of size $k\geq1$. One can use Lagrange inversion formula, with $R(x) = \exp({x\over{1-x}})$.
Alternately, observe that the factor $\frac{(n-1)!}{k!} {n-2\choose k-1}$ in (\ref{eq:ocn})
represents the number of partitions of a set of size $n-1$ into $k$ totally 
ordered parts so that the Pr\"ufer-type bijection of Springer \cite{Sp96} can be used here.
\qed

\begin{proposition}
\emph{\cite{Sp96}}
Let $(n_2,n_3,\ldots)$
be a sequence of non-negative integers and $n=\sum_{i\geq2}n_i(i-1)+1$.
Then the number $\Oc(n_2,n_3,\ldots)$ of oriented cacti on $[n]$ having
$n_i$ cycles of size $i$ for each $i$, is given by
\begin{equation}
\Oc(n_2,n_3,\cdots) = \frac{(n-1)!}{n_2!n_3!\cdots}\, n^{k-1}, \label{eq:ocnn2n3}
\end{equation}
where  $k=\sum_{i\geq 2}n_i$ is the number of cycles.
\end{proposition}

\noindent {\bf Proof.}
Again, it is possible to use Lagrange inversion or the Pr\"ufer-type correspondence
of Springer. However the result now follows simply from equation (\ref{eq:hun2n3bis})
since it is easy to see that
$$
\Oc(n_2,n_3,\ldots) = \prod_{i \geq 2} (i-1)!^{n_i} \Hu(n_2,n_3,\ldots). 
$$
Indeed, a set of size $i$ can be structured into an oriented cycle in $(i-1)!$ ways.
\qed

Finally, let us consider the species $\Ca$ of cacti which is of the form $\CB$, where 
$B=\sum_{k\geq2}P_k$ is the species of polygons. By convention, a polygon of size $2$
is simply an edge ($K_2=E_2$).  See Example 1.1.3 and Figure \ref{fig:cactus}. 
These structures frequently appear in mathematical research, for example, more recently,
in the context of the Traveling Salesman Problem (see for instance \cite{Fl99}) 
and in the classification of graphic matroids \cite {Za02}.

\begin{proposition} \emph{(Ford and Uhlenbeck \cite{FoUhI56})}
Let $(n_2,n_3,\ldots)$
be a sequence of non-negative integers and $n=\sum_{i\geq2}n_i(i-1)+1$.
Then the number $\Ca(n_2,n_3,\ldots)$ of cacti on $[n]$ having
$n_i$ polygons of size $i$ for each $i$, is given by
\begin{equation}
\Ca(n_2,n_3,\cdots) = \frac{1}{2^{\sum_{j\geq3} n_j}}
    \frac{(n-1)!}{\prod_{j\geq2}n_j!}\, n^{k-1}, \label{eq:cann2n3}
\end{equation}
where  $k=\sum_{i\geq 2}n_i$ is the number of polygons.
\end{proposition}

\noindent {\bf Proof.}
Indeed, since any (labelled) polygon of size $\geq3$ has $2$ orientations,
we see that
$$
\Oc(n_2,n_3\ldots) = 2^{\sum_{j\geq3}n_j}\Ca(n_2,n_3,\cdots).
$$
The result then follows from (\ref{eq:ocnn2n3}).
 \qed

As observed in \cite{FoUhI56}, this corrects the formula given in the introduction of \cite{HaUl53}
which is rather Husimi's formula (\ref{eq:hun2n3bis}) for Husimi graphs.
By summing over the polygon-size distribution, we finally obtain: 

\begin{proposition}
The number $\Ca_n=|\Ca[n]|$ of (labelled) cacti on $[n]$, for $n \geq 2$, is given by
\begin{equation}
\Ca_n=  \sum_{k\geq0} \sum_{{n_2+n_3+\cdots=k \atop n_2+2n_3+\cdots=n-1}}  
         \frac{(n-1)! \ n^{k-1}}{2^{n_3+n_4+\cdots} \ n_2!n_3!\cdots} \ . 
\label{eq:can}
\end{equation}
\end{proposition}

\subsection{Unlabelled enumeration}

For unlabelled enumeration, the cases of rooted and of unrooted $\CB$-graphs
must be treated separately. There is a basic species relationship which permits
the expression of the unrooted species in terms of the rooted ones. It plays the role
of the classical Dissimilarity characteristic theorem for graphs (see \cite{HaPa73}).
Let $\B$ be a species of 2-connected graphs and $\CB$, the associated species
of $\CB$-graphs, that is of connected graphs with blocks in $\B$.
We introduce the following notations:
\begin{enumerate}
\item
$\CBB$ is the species of $\CB$-graphs with a distinguished block,
\item
$\CBSB$ is the species of $\CB$-graphs with a distinguished vertex-rooted block.
\end{enumerate}
\begin{theorem}
\emph{Dissymmetry Theorem for Graphs} \emph{\cite{Le88}, \cite{LeMi92}, \cite{BLL98}}. 
The species $\CB$ of connected graphs whose blocks are
in $\B$ and its associated rooted species are related by the following isomophism:
\begin{equation}
\CBpoint + \CBB = \CB + \CBSB. \label{eq:CBB}
\end{equation}
This identity can also be written as 
\begin{equation}
\CBpoint + \B(\CBpoint) = \CB + \CBpoint \cdot \Bprime(\CBpoint). \label{eq:CBBbis}
\end{equation}
\end{theorem} 

\noindent {\bf Proof.}
The proof of (\ref{eq:CBB}) is remarkably simple. It uses the concept of \emph{center}
of a $\CB$-graph, which is defined as the center of the associated block-cutpoint tree
(see Figure \ref{fig:graphblock}, c)). The center will necessarily be either a vertex 
(in fact an articulation point),
or a block of the $\CB$-graph. Now a structure $s$ belonging to the left-hand-side of 
(\ref{eq:CBB}) is a $\CB$-graph which is rooted at either a vertex or a block (a cell). It can 
happen that the rooting is performed right at the center. This is canonically equivalent
to doing nothing and is represented by the first term in the right-hand-side of
(\ref{eq:CBB}). On the other hand, if the rooting is done at an off-center cell, a vertex 
or a block, then there is a unique incident cell of the other kind
(a block for a vertex and vice-versa) which is located towards the center, thus defining
a unique $\CBSB$-structure. It is easily checked that this correspondence is bijective
and independent of any labelling, giving the desired species isomorphism.

 For (\ref{eq:CBBbis}), it suffices to verify the species isomophisms
$\CBB = \B(\CBpoint)$ and $\CBSB = \CBpoint \cdot \Bprime(\CBpoint)$.
Details are left as an exercise. \qed
%

   The method of enumeration of unlabelled $\CB$-graphs then consists in enumerating
first the rooted unlabelled $\CB$-graphs, using the functional equation (\ref{eq:cbpoint})
or (\ref{eq:cbpointw}), and then the unrooted ones, using (\ref{eq:CBBbis}).
Notice that these steps will lead not to explicit but rather to
recursive formulas for the desired numbers or total weights.  Here
we illustrate the method in detail in the case of triangular cacti (Example 1.1.4).
We will also give some results for Husimi graphs and oriented cacti.

\begin{figure}[bht] 
\begin{center}
\includegraphics[height=6cm]{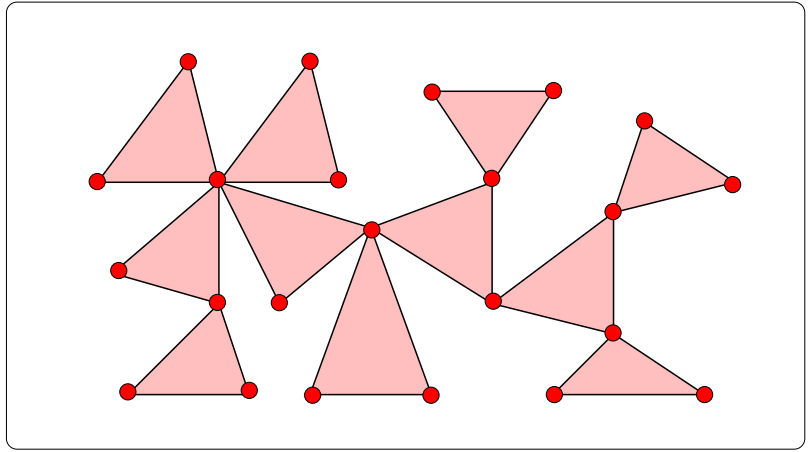}
\end{center}
\caption{Unlabelled triangular cactus} \label{fig:cactustri}
\end{figure}

  Recall that a triangular cactus is a connected graph all of whose blocks are
triangles. See Figure \ref{fig:cactustri} for an example. 
These structures (and also the quadrangular cacti) 
were enumerated by Harary, Norman, and Uhlenbeck \cite{HaNo53}, \cite{HaUl53}. 
Here we go further by giving recurrence formulas for their numbers.  
Let $\delta=\CB$ and $\Delta=\delta^{\bullet}$ denote the species of triangular cacti and
of rooted triangular cacti, respectively. 
Also set 
\begin{equation}
\widetilde{\delta}(x)=\sum_{n\geq1} d_nx^n \ \ 
   \mathop{\rm and} \ \ \widetilde{\Delta}(x)=\sum_{n\geq1}D_nx^n,  \label{eq:deltax}
\end{equation}
where $d_n$ and $D_n$ denote the numbers of unlabelled triangular cacti and
rooted triangular cacti, respectively.
Here we have $\B=K_3=E_3$, and $\Bprime=E_2$.
The functional equations (\ref{eq:cbpoint}) and (\ref{eq:CBBbis}) give
\begin{equation}
\Delta=X \, E(E_{2}(\Delta))  \label{eq:deltapoint}
\end{equation}
and
\begin{equation}
\Delta+E_3(\Delta) = \delta+\Delta \, E_2(\Delta), \label{eq:dissymdelta}
\end{equation}
respectively. Note that 
\begin{eqnarray*}
Z_{E\circ E_2}(x_1,x_2,x_3,\ldots) & = & Z_E(Z_{E_2}(x_1,x_2,\ldots),
                                           Z_{E_2}(x_2,x_4,\ldots),\ldots) \\
 & = & Z_E(\frac{1}{2}(x_1^2+x_2),\frac{1}{2}(x_2^2+x_4),\ldots) \\
 & = & \exp{\left(\sum_{k\geq1} \frac{1}{k}\frac{(x_k^2+x_{2k})}{2}\right)}.
\end{eqnarray*}
From (\ref{eq:deltapoint}), we deduce (see \cite{HaUl53})
\begin{eqnarray}
\widetilde{\Delta}(x) & = & (XE(E_2(\Delta)))^{\sim}(x) \nonumber \\
   & = & x Z_{E\circ E_2}(\widetilde{\Delta}(x),\widetilde{\Delta}(x^2),
              \widetilde{\Delta}(x^3),\ldots) \nonumber \\
   & = & x\exp{\left(\sum_{k\geq1} \frac{1}{2k}\left(\widetilde{\Delta}^2(x^k)
                 + \widetilde{\Delta}(x^{2k})\right)\right)}. \label{eq:deltapointtilde}
\end{eqnarray}
A recurrence formula for $D_n$ can be obtained by taking $x$ $\times$ the derivative of 
(\ref{eq:deltapointtilde}). First set
$$
b(x) = x \frac{d}{dx}\sum_{k\geq1} \frac{1}{2k}\left(\widetilde{\Delta}^2(x^k)
           + \widetilde{\Delta}(x^{2k})\right).
$$
We then have
$x \widetilde{\Delta}^{\prime}(x) = \widetilde{\Delta}(x)+\widetilde{\Delta}(x)b(x)$
and 
\begin{equation}
x\widetilde{\Delta}^{\prime}(x) - \widetilde{\Delta}(x) = \widetilde{\Delta}(x)b(x).
                    \label{eq:xdeltaprime}
\end{equation}
But
\begin{eqnarray*}
b(x) & = & x\frac{d}{dx}\sum_{k\geq1} \frac{1}{2k}\left(\widetilde{\Delta}^2(x^k)
          + \widetilde{\Delta}(x^{2k})\right)\\
     & = & x \cdot \sum_{k\geq1} \frac{1}{2k}\left(2\widetilde{\Delta}(x^k)
            \widetilde{\Delta}^{\prime}(x^k)kx^{k-1}
            +\widetilde{\Delta}^{\prime}(x^{2k})2kx^{2k-1}\right) \\
     & = & \sum_{k\geq1} \left(\widetilde{\Delta}(x^k)\widetilde{\Delta}^{\prime}(x^k)x^{k}
              +\widetilde{\Delta}^{\prime}(x^{2k})x^{2k}\right)\\
     & = & \sum_{k\geq1} \left(\sum_{h\geq1}\sum_{j\geq1} jD_hD_jx^{k(h+j)}
             +\sum_{j\geq1} jD_jx^{2kj}\right)\\
     & = & \sum_{m\geq1}\left( \sum_{{d|m \atop d\geq2}} \sum_{j=1}^{d} jD_{d-j}D_j 
             + \sum_{{d|m \atop d \ \mbox{\tiny{even}}}} \frac{d}{2}D_{\frac{d}{2}}\right) x^m.
\end{eqnarray*}
By extracting the coefficient of $x^{n+1}$ in (\ref{eq:xdeltaprime}), we obtain

\begin{proposition}
The numbers $D_n$, of unlabelled rooted triagular cacti on $n$ vertices satisfy
$D_1=1$ and the recurrence formula, for $n\geq1$,
\begin{equation}
D_{n+1}=\frac{1}{n}\sum_{m=1}^{n}D_{n-m+1}\left( \sum_{d|m \atop d\geq2} \sum_{j=1}^{d} jD_{d-j}D_j 
                + \sum_{{d|m \atop d \ \mbox{\tiny{even}}}} \frac{d}{2}D_{\frac{d}{2}}\right).
                    \label{eq:recurdn}
\end{equation}
\end{proposition}

In order to enumerate unlabelled \emph{unrooted} triangular cacti, we use (\ref{eq:CBBbis}),
that is,
$$
\delta=\Delta+E_3(\Delta)-\Delta E_2(\Delta).
$$
The passage to the tilde generating functions yields (see \cite{HaNo53})
\begin{eqnarray}
\widetilde{\delta}(x) & = & \widetilde{\Delta}(x)+(E_3(\Delta(x)))^{\sim} 
           -  \widetilde{\Delta}(x) (E_2(\Delta(x)))^{\sim} \nonumber \\
 & = & \widetilde{\Delta}(x) + Z_{E_3}(\widetilde{\Delta}(x),\widetilde{\Delta}(x^2),\ldots) 
    - \widetilde{\Delta}(x)Z_{E_2}(\widetilde{\Delta}(x),\widetilde{\Delta}(x^2),\ldots) \nonumber\\ 
 & = & \widetilde{\Delta}(x)+ \frac{1}{6}(\widetilde{\Delta}^{3}(x) 
        + 3\widetilde{\Delta}(x)\widetilde{\Delta}(x^2)+ 2\widetilde{\Delta}(x^3)) 
        -  \frac{1}{2}\widetilde{\Delta}(x)(\widetilde{\Delta}^2(x) 
                     +\widetilde{\Delta}(x^2)) \nonumber \\
 & = & \widetilde{\Delta}(x) + \frac{1}{3}\left(\widetilde{\Delta}(x^3)
                     -\widetilde{\Delta}^3(x)\right).  \label{eq:dissymdeltax}
\end{eqnarray}
By extracting the coefficient of $x^{n}$, we finally obtain
%

\begin{proposition}
The numbers $d_n$, of unlabelled (unrooted) triangular cacti on $n$ vertices satisfy, 
for $n\geq1$,
\begin{equation}
d_n=D_n +\frac{1}{3}\left(\mbox{\Large{$\chi$}}(3|n) D_{\frac{n}{3}}
              -\sum_{i+j+k=n} D_iD_jD_k\right).   \label{eq:dn}
\end{equation}
\end{proposition}

   We now turn our attention to the species $\Hu$ of Husimi graphs. Let us denote by
$h_n$ and $H_n$, the numbers of unlabelled Husimi graphs and rooted Husimi graphs,
respectively, and set
\begin{equation}
h(x) = \sum_{n\geq1}h_nx^n = \widetilde{\Hu}(x) \ \ \mathop{\rm and} \ \ 
H(x) = \sum_{n\geq1}H_nx^n = \widetilde{\Hup}(x).   \label{eq:hux}
\end{equation}
Here, as seen in Section 3.1,
$\B = E_{\geq 2}$ and $\B^{\prime} = E_{\geq 1} = E_+$, and the basic functional equation
(\ref{eq:hupoint}) translates, for the tilde generating function, into (see \cite{No54}, p. 51)
\begin{equation}
H(x) =x \exp{\sum_{k\geq1} \frac{1}{k}\left(\exp
         {\left(\sum_{m\geq1} \frac{1}{m} H(x^{mk})\right)-1}\right)}.   \label{eq:huptilde1}
\end{equation}
We introduce the auxiliary series
\begin{equation}
b(x) =  \sum_{n\geq1}b_nx^n = x\frac{d}{dx}\sum_{k\geq1} 
    \frac{1}{k}\left(\exp \left(\sum_{m\geq1} \frac{1}{m} H(x^{mk}) \right)-1\right), 
   \label{eq:bhux}
\end{equation}
\begin{equation}
\varphi(x) = \sum_{n\geq1} \varphi_nx^n = 
               x\exp{\left(\sum_{j\geq1} \frac{1}{j}\Hn(x^j)\right)}.
   \label{eq:phix}
\end{equation}
Then, after some computations, similar to those for triangular cacti, we find:

\begin{proposition}
The numbers $H_n$, of unlabelled rooted Husimi graphs on $n$ vertices can be computed
by the following recusive scheme: \hspace{.1em} $\varphi_1 = 1$, \hspace{.1em} $H_1 = 1$,
\hspace{.1em} and, for $n\geq1$,
\begin{eqnarray}
b_n & = & \sum_{d|n} \sum_{h=1}^{d} \sum_{l|d-h+1} lH_l\ \varphi_h,
                             \label{eq:bnreponse} \\
H_{n+1} & = & \frac{1}{n} \sum_{k=1}^{n} H_{n-k+1} b_k,
                             \label{eq:recurrenceHustilde} \\
\varphi_{n+1} & = & \frac{1}{n}\sum_{m=1}^{n} \varphi_{n-m+1} \sum_{d|m} dH_d.
                             \label{eq:varphin} 
\end{eqnarray}
\end{proposition}

For unlabelled \emph{unrooted} Husimi graphs, we use the Dissymmetry Theorem, in the form 
(\ref{eq:CBBbis}), which gives
\begin{eqnarray*}
\Hup + E_{\geq2}(\Hup) & = & \Hu + \Hup\cdot E_+(\Hup), \\
E_ + (\Hup) & = & \Hu+\Hup\cdot E_+(\Hup), 
\end{eqnarray*}
and finally

\begin{equation}
\Hu = (E_+(\Hup)) \cdot (1-\Hup).  \label{eq:relationHus2}
\end{equation}

For the tilde generating series, we deduce that
\begin{eqnarray*} 
h(x) & = &  \left(\exp{\left(\sum_{k\geq1} \frac{1}{k}H(x^k)\right)}-1\right)
                   \cdot \left(1-H(x)\right)  \\
 & = &  \left(\frac{\varphi(x)}{x}-1\right)\cdot\left(1-H(x)\right) , 
\end{eqnarray*}
and we obtain the following result.

\begin{proposition}
The numbers $h_n$, of unlabelled Husimi graphs on $n$ vertices satisfy
\begin{equation}
h_n = \varphi_{n+1} - \sum_{k=1}^{n-1} \varphi_{k+1} H_{n-k} \label{eq:hn},
\end{equation}
where the numbers $H_n$ and $\varphi_n$ are given by Proposition 3.10.
\end{proposition}

The same method can be applied to weighted $\CB$-graphs where the weight is defined
by (\ref{eq:yiweight}), that is, with the variables $(y_2,y_3,\ldots)$ marking the block sizes.
We have done the computations for the weighted species $\Oc_w$, of oriented cacti,
where $\B_w = \sum_{m\geq2}(C_m)_{y_m}$ and $\Bprime_w = \sum_{m\geq1}(L_m)_{y_{m+1}}$,
where $(C_m)_{y_m}$ denotes the species of oriented cycles of length $m$ and weight $y_m$,
and $L_m=X^m$ is the species of lists of size $m$, endowed here with the weight $y_{m+1}$.
Introduce the following tilde generating series, whose coefficients are polynomials
in the variables $\yg = (y_2,y_3,\ldots)$:
\begin{eqnarray*}
o(x;\yg) := \widetilde{\Oc_w}(x) & = & \sum_{n\geq1} o_n({\yg})x^n, \\
O(x;\yg) := \widetilde{\Ocp_w}(x) & = & \sum_{n\geq1} O_n({\yg})x^n. 
\end{eqnarray*}
The two basic species functional equations are
\begin{equation}
\Ocp_w =  X\cdot E\left(\sum_{m\geq1} (L_m)_{y_{m+1}}(\Ocp_w)\right), \label{eq:ocpw}
\end{equation}
\begin{equation}
C_w(\Ocp) = \Oc + L_{\geq2,w}(\Ocp), \label{eq:ocwdiss}
\end{equation}
where $C_w = \sum_{m\geq1} (C_m)_{y_m}$, with $C_1=X, \ y_1=1$, and $L_{\geq2,w} = \sum_{m\geq2}(X^m)_{y_m}$.
%
Note that 
$$
 Z_{C_w} = \sum_{m\geq1} \frac{y_m}{m} \sum_{d|m} \phi(d)x_d^{\frac{m}{d}}
         = \sum_{d\geq1} \sum_{h\geq1} \frac{y_{dh}}{dh}\phi(d)x_d^h, 
$$
where $\phi$ is Euler's totient function, and that $Z_{L_{\geq2,w}} = \sum_{m\geq2}y_mx_1^m$.
From (\ref{eq:ocpw}), we deduce, using the plethystic composition rule for the tilde generating function
of weighted species (see (4.3.1) of \cite{BLL98}),
\begin{equation}   
O(x;\yg)=x\exp{\left(\sum_{k\geq1} \frac{1}{k}\sum_{j\geq1} y_{j+1}^k O^j(x^k;\yg^k)\right)}, 
 \label{eq:ocpwtilde}  
\end{equation}
where $\yg^k=(y_2^k,y_3^k,\ldots)$. We then set
\begin{equation}   
O^j(x;\yg) = (O(x;\yg))^j = \sum_{n\geq1} O_n^{(j)}({\yg})x^n, 
\end{equation}
\begin{equation}   
b(x;\yg) = \sum_{m\geq1} b_m(\yg) x^m = x\parx \sum_{k\geq1} \frac{1}{k}\sum_{j\geq1} y_{j+1}^{k}
                 O^j(x^k;\yg^k).  \label{eq:bxy}  
\end{equation}
From (\ref{eq:ocwdiss}), we deduce 
\begin{equation}  
o(x;\yg) = \sum_{d\geq1} \phi(d) \sum_{h\geq1} \frac{y_{dh}}{dh} O(x^d;\yg^d)
           - \sum_{m\geq2}y_m O(x;\yg)^m
\label{eq:ocwdissx}
\end{equation}
and we obtain, after some computations, the following result.

\begin{proposition}
The generating polynomials $O_n(\yg)$ and $o_n(\yg)$ for 
unlabelled rooted and unrooted (resp.) oriented cacti on $n$ vertices are given 
by the following recursive scheme:

\begin{equation}   
O_{n+1}(\yg) = \frac{1}{n} \sum_{m=1}^n O_{n-m+1}(\yg) b_m(\yg), \label{eq:recuropny}  
\end{equation}
where
\begin{equation}   
b_m(\yg) = \sum_{d|m} \sum_{i=1}^{d} \sum_{j=1}^{d-i+1} 
ij\ y_{j+1}^{\frac{m}{d}}O_i(\yg^{\frac{m}{d}})O_{d-i}^{(j-1)}(\yg^{\frac{m}{d}}),  \label{eq:bmy}
\end{equation}
and
\begin{equation}   
o_n(\yg) = \sum_{d|n} \phi(d) \sum_{h=1}^{\frac{n}{d}} \frac{y_{dh}}{dh} 
             O_{\frac{n}{d}}^{(h)}(\yg^d) - \sum_{m=2}^n y_m O_n^{(m)}(\yg).  \label{eq:ony}
\end{equation}
\end{proposition}

\section{Molecular expansions}

The molecular expansion of a species $F$ is a description and a classification of the 
unlabelled $F$-structures according to their stabilizers, within the language of species.
It is interesting and useful to have at hand the first few terms of the molecular
expansion of $\CB$-graphs. This is possible with the Maple package Devmol,
developped at LaCIM. See \cite{AuLe03}.

As an example, consider again the species $\Hu_w$ of Husimi graphs, weighted
according to their block-size distribution. Here $\B_w = E_{\ge 2,w} = \sum_{k \ge 2} (E_k)_{y_{k}}$.
This is the species $BB$ in the following extract from a Maple session, using the package Devmol.
The Devmol procedure "\texttt{CBgraphes}(BB,n)" produces the molecular 
expansion of the species of $\Gcon_{BB}$-graphs, up to the given truncation order $n$. 
The expansion is then collected by degree (vertex number).
The weight variables $y_2, y_3, \ldots, y_n$ are
first declared to Devmol. Afterwards, the monomial weights in the $y_i$'s appear multiplicatively
in the expressions. Here is the extract:

\begin{maplelike}
\mwsin{n := 6;}
\mwsout{%
$$n := 6$$
}
\mwsin{ajoutvv(seq(y[k], k = 1..n));}
\mwsout{%
$$\{t, \,{y_{1}}, \,{y_{2}}, \,{y_{3}}, \,{y_{4}}, \,{y_{5}}, \,{y_{6}}\}$$
}
\mwsin{BB := sum(y[k]*E[k](X), k = 2..n);}
\mwsout{%
$$\begin{array}{c}
BB := {y_{2}}\,{E_{2}}(X) + {y_{3}}\,{E_{3}}(X) + {y_{4}}\,{E_{4}}
(X) + {y_{5}}\,{E_{5}}(X) + {y_{6}}\,{E_{6}}(X)
\end{array}$$
}
\mwsin{HuW := CBgraphes(BB,n):}
\mwsin{affichertable(tablephom(HuW));}
\mwsout{%
$$\begin{array}{c}
1\\
X
\end{array}$$
}
\mwsout{%
$$\begin{array}{c}
2\\
{y_{2}}\,{E_{2}}(X)
\end{array}$$
}
\mwsout{%
$$\begin{array}{c}
3\\
X\,{y_{2}}^{2}\,{E_{2}}(X) + {y_{3}}\,{E_{3}}(X)
\end{array}$$
}
\mwsout{%
$$\begin{array}{c}
4\\
{y_{2}}\,X^{2}\,{y_{3}}\,{E_{2}}(X) + X\,{y_{2}}^{3}\,{E_{3}}(X)
 + {y_{2}}^{3}\,{E_{2}}(X^{2}) + {y_{4}}\,{E_{4}}(X)
\end{array}$$
}
\mwsout{%
$$\begin{array}{c}
5\\
({y_{3}}\,{y_{2}}^{2} + {y_{2}}^{4})\,X^{3}\,{E_{2}}(X) + X\,{y_{2}}^{2}\,{E_{2}}(X)^{2}\,{y_{3}} 
+ ({y_{3}}\,{y_{2}}^{2} + {y_{2}}^{4})\,X\,{E_{2}}(X^{2}) \\
+ X\,{y_{2}}^{4}\,{E_{4}}(X)  + X\,{y_{3}}^{2}\,{E_{2}}({E_{2}}(X)) + {y_{5}}\,{E_{5}}(X) 
+ {y_{2}}\,X^{2}\,{y_{4}}\,{E_{3}}(X)

\end{array}$$
}
\mwsout{%
$$\begin{array}{c}
6\\
{y_{2}}^{5}\,X^{2}\,{E_{2}}(X^{2}) + ({y_{2}}^{2}\,{y_{4}} + {y_{3}}\,{y_{4}} + {y_{3}}\,{y_{2}}^{3})\,X\,{E_{2}}(X)\,{E_{3}}(X) \\
 + {y_{4}}\,{E_{2}}(X)\,{y_{2}}^{2}\,{E_{2}}(X^{2})
\mbox{} + ({y_{3}}^{2}\,{y_{2}} + {y_{2}}^{5} + 3\,{y_{3}}\,{y_{2}}^{3})\,X^{4}\,{E_{2}}(X) \\
+ {y_{2}}^{5}\,{E_{2}}(X^{3}) + ({y_{3}}^{2}\,{y_{2}} + {y_{2}}^{5})\,{E_{2}}(X\,{E_{2}}(X))
\mbox{} + {y_{2}}\,{y_{3}}^{2}\,X^{2}\,{E_{2}}({E_{2}}(X)) \\
+ {y_{3}}\,{y_{2}}^{3}\,{E_{3}}(X^{2}) + {y_{2}}\,{y_{5}}\,X^{2}\,{E_{4}}(X)
 + {y_{6}}\,{E_{6}}(X) + {y_{3}}\,{y_{2}}^{3}\,X^{2}\,{E_{2}}(X)^{2} \\
\mbox{} + {y_{3}}\,{y_{2}}^{3}\,X^{6} + ({y_{2}}^{5} + {y_{2}}^{2
}\,{y_{4}})\,X^{3}\,{E_{3}}(X) + {y_{2}}^{5}\,X\,{E_{5}}(X) 
\end{array}$$
}
\end{maplelike}

\vspace{.5\baselineskip}
\begin{figure}[th]
\centerline{\input{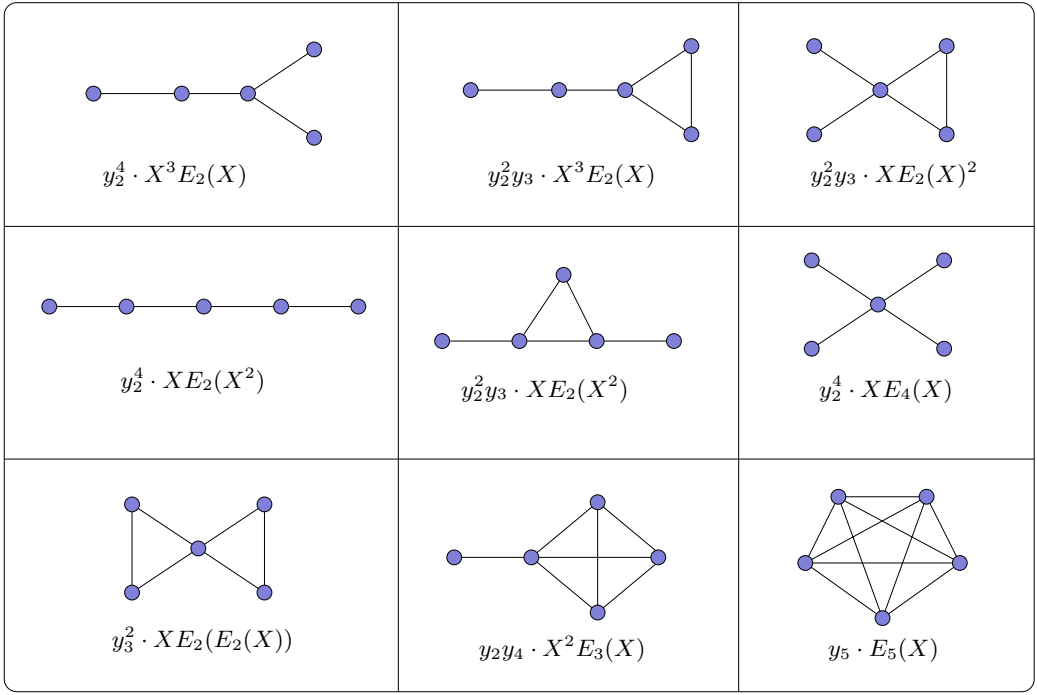}}
  \caption{Husimi graphs of size $5$}
  \label{fig:husimi5}
\end{figure}
\begin{figure}[th]
\includegraphics{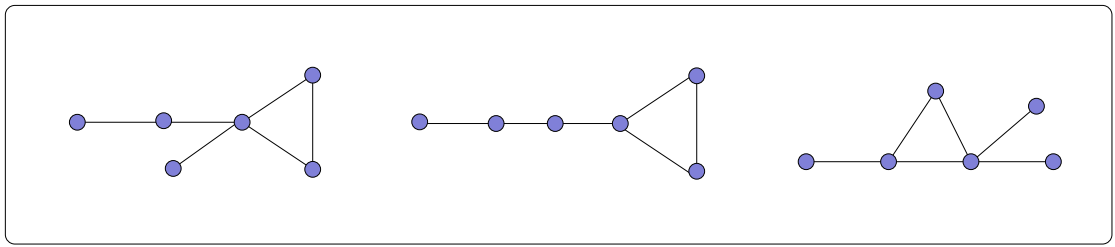}
  \caption{The three unlabelled Husimi graphs of molecular type $y_2^3y_3\cdot X^4E_2(X)$}
  \label{fig:husimi6}
\end{figure}

As an illustration, the nine terms of size $5$ are shown in Figure \ref{fig:husimi5}.
Observe that multiciplicities will frequently occur in these expansions.
For size $6$, for example, the term $3\,y_3\,y_2^3\,X^4E_2(X)$ corresponds 
to the $3$ similar types of Husimi graphs shown in Figure \ref{fig:husimi6}.
{\small 

}
\noindent Pierre Leroux \\
LaCIM, D\'ep. de math\'ematiques,\\ Universit\'e du Qu\'ebec \`a Montr\'eal,\\
C.P. 8888, Succ. Centre-Ville,\\Montr\'eal (Qu\'ebec), Canada  H3C 3P8.\\
leroux.pierre@uqam.ca, www.lacim.uqam.ca/$\sim$leroux

\end{document}